\documentclass{amsart}
\usepackage{amsfonts,amssymb,amscd,amsmath,enumerate,verbatim,calc}
\usepackage[all]{xy}
\usepackage{hyperref}

\newcommand{\CM}{CM}
\newcommand{\MCM}{MCM}
\newcommand{\IFF}{\text{if and only if}}
\newcommand{\wrt}{with respect to}
\newcommand{\sg}{s.g}

\newcommand{\Z}{\mathbb{Z} }
\newcommand{\n}{\mathfrak{n} }
\newcommand{\m}{\mathfrak{m} }
\newcommand{\q}{\mathfrak{q} }

\newcommand{\xar}{\longrightarrow}
\newcommand{\ov}{\overline}

 \newcommand{\rt}{\rightarrow}
\newcommand{\ve}{\varepsilon}

\newcommand{\om}{\omega}

\newcommand{\y}{\mathbf{y} }
\newcommand{\z}{\mathbf{z} }
\newcommand{\trm}{\mathbf{m} }
\newcommand{\xbold}{\mathbf{x} }
\newcommand{\C}{\mathbf{C}_{\bullet} }

\newcommand{\K}{\mathbf{K}_{\bullet}}

\newcommand{\post}{\operatorname{post}}

\newcommand{\pd}{\operatorname{projdim}}

\newcommand{\type}{\operatorname{type}}
\newcommand{\Ass}{\operatorname{Ass}}

\newcommand{\projdim}{\operatorname{projdim}}

\newcommand{\ann}{\operatorname{ann}}
\newcommand{\grade}{\operatorname{grade}}
\newcommand{\depth}{\operatorname{depth}}
\newcommand{\image}{\operatorname{image}}
\newcommand{\socle}{\operatorname{socle}}
\newcommand{\Syz}{\operatorname{Syz}}
\newcommand{\Tor}{\operatorname{Tor}}
\newcommand{\Hom}{\operatorname{Hom}}
\newcommand{\Ext}{\operatorname{Ext}}

\theoremstyle{plain}
\newtheorem{thm}{Theorem}

\newtheorem{theorem}{Theorem}[section]
\newtheorem{corollary}[theorem]{Corollary}
\newtheorem{lemma}[theorem]{Lemma}
\newtheorem{proposition}[theorem]{Proposition}

\theoremstyle{definition}
\newtheorem{definition}[theorem]{Definition}
\newtheorem{s}[theorem]{}
\newtheorem{remark}[theorem]{Remark}
\newtheorem{example}[theorem]{Example}
\newtheorem{construction}[theorem]{Construction}

\begin{document}

\title[ Hilbert Function]{ The Hilbert Function of a Maximal Cohen-Macaulay
Module}
\author{Tony ~J.~Puthenpurakal}

 \begin{abstract}
We study Hilbert functions of maximal \CM \ modules over \CM \ local rings.
 We show that if $A$ is a hypersurface ring with dimension $d > 0$ then the
Hilbert function of $M$ \wrt \ $\m$ is non-decreasing. If $A =
Q/(f)$ for some regular local ring $Q$, we determine a lower bound
for $e_0(M)$ and $e_1(M)$ and analyze the case when equality holds.
 When $A$ is Gorenstein  a  relation
between the second Hilbert coefficient of $M$, $A$ and $S^A(M)=  (\Syz^{A}_{1}(M^*))^*  $ is
found when $G(M)$ is \CM \ and $\depth G(A) \geq d-1$. We give
bounds for the first Hilbert coefficients of the canonical module
of a \CM \  local ring and analyze when equality holds. We also give
 good bounds on Hilbert coefficients of $M$ when $M$ is maximal
\CM \ and $G(M)$ is \CM.
\end{abstract}
\date{\today}
 \maketitle

\section*{Introduction}
Let $(A,\m)$ be a $d$-dimensional Noetherian local ring  and $M$ a
finite $A$-module. Let $G(A) = \bigoplus_{n\geq 0}\m^n/\m^{n+1}$
be the associated graded module of $A$ and $G(M) =  \bigoplus
_{n\geq 0}\m^nM /\m^{n+1}M$    the associated graded module of $M$
considered as a $G(A)$-module. We set $ \depth G(M) =
\grade(\mathcal{M},G(M)) $   where $\mathcal{M} = \bigoplus_{n\geq 1}\m^n/\m^{n+1} $
  is the irrelevant maximal ideal of $G(A)$.
If $N$ is an $A$- module then $\mu(N)$ denotes its minimal number
of generators and $\lambda(N)$ denotes its length. The Hilbert
function of $M$ (\wrt \ $\m$) is the function
\[
  H(M,n) = \lambda(\m^{n}M/\m^{n+1}M) \text{ for all } n \geq 0
\]

In this paper we study Hilbert functions of maximal \CM \ (= MCM) modules.
If $A$ is regular then all  \MCM \ modules are free. The next
case is that of a hypersurface ring.

\begin{thm}
   \label{c1thMo}
     Let $(A,\m)$ be a hypersurface ring of positive dimension. If $M  $  is
a \MCM \  $A$-module, then
     the Hilbert function of $ M  $ is non-decreasing.
    \end{thm}
 
This result is a corollary of a more general result (see Theorem \ref{thMo}) which also implies that the
Hilbert function of a complete intersection of codimension 2 and positive dimension is non-decreasing
(see Corollary \ref{c2thMo}). 
 Another application of  Theorem \ref{thMo} yields that, if  $(A,\m)$ is equicharacteristic local ring of
 dimension $d > 0$, $I$ is an $\m$-primary ideal with $\mu(I) = d+ 1 $ and  
$M$ is an \MCM \ $A$-module, then the Hilbert function of $M$ \wrt \ $I$ is non-decreasing
(see Theorem \ref{muI}).

    The  formal power series
\[
  H_M(z) = \sum_{n \geq 0} H(M,n)z^n
\]
is called the \emph{Hilbert series } of $M$. It is well known that
it is of the form
\begin{equation*}
 H_M(z)  = \frac{h_M(z)}{(1-z)^r}, \ \text{where}\ \ r = \dim M \ \  \
 \text{and} \ 
 h_M(z)  \in
\mathbb{Z}[z].
\end{equation*}
We call  $h_M(z)$ the \emph{h-polynomial} of $M$. 
 If $f$ is a polynomial we use $f^{(i)}$ to denote its  $ i $-th derivative.
 The integers $ e_i(M) = h^{(i)}_{M} (1)/i!$ for $ i \geq 0 $  are called
the
 \emph{Hilbert coefficients} of $M$. The number $ e(M) = e_0(M) $  is
the \emph{multiplicity}
 of $ M $. Set 
$$\chi_i(M) = \sum_{j =0}^{i}(-1)^{i-j} e_{i-j}(M)  + (-1)^{i+1}\mu(M) \quad \text{for each  $i \geq 0$}.$$ 

Let $M$ be a  \MCM \
module  over a hypersurface ring $A = Q/(f)$, where $(Q,\n)$ is a regular local.
If  $0 \rt Q^{n} \xrightarrow{\phi_M} Q^{n} \rt M \rt 0 $  is a minimal presentation of
$M$
then  $ i(M) = \max\{i \mid \ \text{all entries of } \ \phi \ \text{are in} \  \n^i \}$ is an  invariant of $M$.

\begin{thm}
\label{mtHy}
Let $(Q,\n)$ be a regular local ring, $f \in \n^e \setminus\n^{e+1}$, $e \geq 2$, $A= Q/(f)$,
  $M$ a \MCM \ $A$-module
and $K = \Syz^{A}_{1}(M)$. Then
\begin{enumerate}[\rm \quad 1.]
\item
     $e(M) \geq \mu(M)i(M)$ and $e_1(M) \geq \mu(M)\binom{i(M)}{2}$.
\item
 $M$ is a free $A$-module \IFF \ $i(M) =e$.
\item
If $i(M) = e-1$ then $G(M)$ is \CM.
\item
  The following conditions are equivalent:
  \begin{enumerate} [\rm \quad  i.]
  \item
  $e(M) = \mu(M)i(M)$.
 \item
   $e_1(M) = \mu(M)\binom{i(M)}{2}$.
    \item
     $G(M)$ is \CM \ and $h_M(z) = \mu(M)(1+z+\ldots+ z^{i(M) -1})$.
\end{enumerate}
If these conditions hold and $M$ is not free, then $G(K)$ is
 \CM \  and
 $h_K(z) = \mu(M) (1+z+\ldots+ z^{e - i(M) -1}).$
\end{enumerate}
\end{thm}

If $A$ is complete and $M$ is a Cohen-Macaulay(= \CM) $A$-module then there exists a Gorenstein local ring $R$
such that $M$ is a \MCM \ $R$-module.  So it is  significant to see how Hilbert functions
of \MCM \ modules over Gorenstein rings behave. 

  When $A$ is \CM, $M$ is a \MCM \ $A$-module and $N =  \Syz^{A}_{1}(M)$ then
  \begin{equation}
\label{cricket1}
  \mu(M) e_1(A) \geq  e_1(M) + e_1(N) \quad \& \quad  \mu(M) \chi_1(A) \geq \chi_1(M) + \chi_1(N).
  \end{equation}
For the first inequality see \cite[17(3)]{Pu1}. The second follows from  \cite[21(1)]{Pu1}.

   In the theorem below we  establish similar inequalities for higher
Hilbert coefficients of \MCM \ modules over Gorenstein
rings. For every $A$-module we set $M^* = \Hom_A(M,A)$. Note that
if $M$ is \MCM \ then so is $M^*$ cf.\cite[3.3.10.d]{BH}.
Also, $\type(M) = \dim_k \Ext_{A}^{d}(k,M) $ denotes the Cohen-Macaulay
type of $M$.
 \begin{thm}
 \label{the2}
  Let $(A,\m)$ be a Gorenstein local ring. Let   $M  $  be a \MCM \
$ A $-module. Set $\tau = \type(M)$ and
   $S^A(M) = \left(\Syz^{A}_{1}(M^*)\right)^*$. If $ G(M) $  is \CM \ and
$\depth G(A) \geq d - 1$ then the following
   hold
   \begin{enumerate}[\quad \rm 1.]
   \item
    $\tau e_2(A) \geq  e_2(M) + e_2(S^{A}(M)) $ and $\tau \chi_2(A) \geq  \chi_2(M) + \chi_2(S^{A}(M)) $.
    \item
     $\type(M)e_i(A) \geq  e_i(M)$   and $\type(M)\chi_i(A) \geq  \chi_i(M)$  for each $i\geq 0$.
  \end{enumerate}
 \end{thm}

Let $A$ be \CM \ with 
a canonical module $\om_A$. 
 Set $\tau = \type A$.
 It is well known that
  $e_{0}(\om_A) = e_{0}(A)$.
Using  \cite[Theorem 18]{Pu1} it follows that $  e_1(\om_A) \leq \tau e_1(A)$ with equality
\IFF \ $A$ is Gorenstein. Here we give a lower bound on $e_1(\om_A)$. 

\begin{thm}
\label{omega}
Let $(A,\m)$ be a \CM \ local ring of dimension $d \geq 1$ and with a canonical module $\om_A$. Set $\tau = \type A$. We have
\begin{enumerate}[\rm (1.)]
\item
$ \tau^{-1}e_1(A) \leq e_1(\om_A) \leq \tau e_1(A).$
\item
\begin{enumerate}[\rm (a.)]
\item
$e_1(\om_A) = \tau e_1(A)$ iff $A$ is Gorenstein.
\item
$e_1(A) = \tau e_1(\om_A)$ iff $A$ is Gorenstein or $A$ has minimal multiplicity.
\end{enumerate}
\item
If $\dim A = 1 $ and $G(A)$ is \CM \ then

\noindent$\displaystyle{ e_i(A) \leq \tau e_i(\om_A) \quad \text{and} \quad \chi_i(A) \leq \tau \chi_{i}(\om_A)  \text{ for each} \  i \geq 1.} $
\end{enumerate}
\end{thm}

Let $A = k[[x_1,\ldots,x_n]]/\q$ be \CM\  with $\q \subseteq (\xbold)^2$ and
$k$ an infinite field. Let $1\leq r \leq d $.  For any two sets of $ r$ sufficiently general $k$-linear
combinations of $ x_1, \ldots, x_\nu $  say $y_1, \ldots, y_r$ and $
z_1, \ldots, z_r $ we show 
$ H(A/(\y) ,n) = H(A/(\z),n)$ for each $n \geq 0$ (see \ref{lem:suffgen}).
 We use it to bound Hilbert coefficients of 
$M$ if $G(M)$ is \CM \  (see  Theorem \ref{existR}).

 Here is an overview of the contents of the paper. 
In Section 1 we introduce notation and discuss a few preliminary facts 
that we need. The proof of the Theorems 1 and 3   involves a study
 of the modules;
  $L_t(M) = \bigoplus_{n\geq 0}  \Tor^{A}_{t}(A/\m^{n+1}, M)$  for all
 $t \geq 0$.
  If $x_1,\ldots,x_s$  is a sequence of elements in  $\m$ then, in section 2, we give
$L_t(M)$ a  structure of a graded
 $A[X_1,\ldots,X_s]$-module. 
 We prove Theorem 1 in Section 3, Theorem 2 in Section 4 and Theorem
3 in Section 5. We prove Theorem 4 in section 6. In Section 7 we prove 
Lemma \ref{lem:suffgen} and use  to prove  Theorem \ref{existR}.

\section{Preliminaries}
 In this paper all rings are Noetherian and all modules are assumed finite 
i.e., finitely generated. 
Let $(A,\m)$ be
 a local ring of dimension $d$ with residue field $k = A/\m$. Let $M$ be
 an
$A$-module. If $m$ is a non-zero
 element of $M$ and if $j$ is the largest integer such that $m \in \m^{j}M$,
then we let $m^*$ denote the image of $m$
 in $\m^{j}M/\m^{j+1}M$. If $L$ is a submodule of $M$, then $L^*$ denotes
the graded submodule of $G(M)$ generated by all
 $l^*$ with $l \in L$. It is well known that $G(M)/L^* = G(M/L)$. An element
$x \in \m$ is said to be \emph{superficial}
  for $M$ if there exists an integer $c > 0$ such that
$$ (\m^{n}M \colon_M x)\cap\m^cM =  \m^{n-1}M \ \text{ for all }\quad n > c.
$$
 Superficial elements always
exist if  $k$ is infinite  \cite[p.\ 7]{Sa}. A
sequence $x_1,x_2,\ldots,x_r$ in a local ring $(A,\m)$ is said to
be a \emph{superficial sequence} for $M$ if $x_1$ is superficial
for $M$ and $x_i$ is superficial for $M/(x_1,\ldots,x_{i-1})M$ for
$2\leq i \leq r$.
 
 \begin{remark}
\label{Pr0}
 If the residue field of $A$ is finite then we resort to the
standard trick to replace $A$ by $A' = A[X]_S$ and $M$ by $M'= M
\otimes_A A'$  where $S =  A[X]\setminus \m A[X]$. The residue
field of $A'$ is $k(X)$, the field of rational functions over $k$.
Furthermore
\begin{equation*}
H(M',n) = H(M,n) \ \  \forall n \geq 0 \ \  \text{and} \  \   \depth_{G(A')}G(M') = \depth_{G(A)}G(M). 
 \end{equation*}
 Clearly $\pd_{A'}M' = \pd_A M $ If $A$ is a  Gorenstein (hypersurface) ring then $A'$ is also Gorenstein (hypersurface) ring.
If $A$ has a canonical module $\om_A$ then 
 $A'$ also  has a canonical module $\om_{A'} \cong \om_A \otimes A'$;  cf. \cite[Theorem 3.3.14]{BH}.  
\end{remark}

Below we collect some basic results needed in the paper. For proofs
 see \cite{Pu1}.
\begin{remark}
 \label{Pr2}
   Let $(A,\m)$ be a Cohen-Macaulay local ring with $\dim A = d > 0$. Let
$M$ be a finite \CM \ $A$-module of dimension $r$.
   Let $x_1,\ldots,x_s$ be a superficial sequence in $M$ with $s \leq r$ and
set $J= (x_1,\ldots,x_s)$.
   The local ring $(B,\n) = (A/J,\m/J)$ and $B$-module $N =
M/JM$  satisfy:
   \begin{enumerate} [\quad \rm 1.]
\item
 $x_1,\ldots,x_s$ is a $M$-regular sequence in $A$.
 \item
 $N$ is a \CM \ $B$-module.
 \item
 $e_i(N) = e_i(M)$ for $i = 0,1,\ldots,r-s$.
 \item
 When $s=1$,  set $x =x_1$ and $ b_n(x, M) = \lambda(\m^{n+1}M \colon_M x)/\m^nM) $. We have:
 \begin{enumerate}[\rm a.]
\item
$b_0(x,M) = 0$ and $ b_n(x,M) = 0 $ for all $ n\gg 0 $.
\item
$ H(M,n) = \sum_{i=0}^{n}H(N,i)  - b_n(x,M)$
\item
$e_r(M) = e_r(N) - (-1)^r\sum_{i \geq 0} b_{i}(x,M)$.
\item
$x^*$ is $G(M)$-regular \IFF \ $b_n(x,M) = 0$ for all $n \geq 0$.
\end{enumerate}
\item
\begin{enumerate}[\rm a.]
\item
$\depth G(M) \geq s$ \IFF \ $x_1^*,\ldots,x_s^*$ is a
$G(M)$ regular sequence.
\item
(\emph{Sally descent}) $\depth G(M) \geq s +1$ \IFF \
$\depth
G(N) \geq 1$.
\end{enumerate}
 \item
 If $\dim M = 1$ then
 set $\rho_n(M) = \lambda(\m^{n+1}M/x\m^nM)$.   We have
   \begin{enumerate} [\rm a.]
    \item
    $ H(M,n) = e(M) - \rho_n(M)$.
    \item
    $ e_i(M) = \sum_{j\geq i -1}\binom{j}{i-1}\rho_j(M) \geq 0 $ for all $i \geq 1$.
   \end{enumerate}
 \item
 If   $x_1,\ldots,x_s$ is also $A$-regular then $\Syz^{B}_{1}(N)
\cong \Syz^{A}_{1}(M)/J \Syz^{A}_{1}(M)$
\item
$\depth G(M) \geq s $ \IFF \ $h_M(z) = h_N(z)$.
\item
$M$ has minimal multiplicity \IFF \ $\chi_1(M) = 0$.
\end{enumerate}
 \end{remark}

\begin{remark}
If    $ \phi : (A,\m) \xar (B,\n)$  is a surjective map of local rings and if $ M $ is a finite $  B$-module then
$ \m^nM = \n^nM  $  for all $ n \geq 0 $. Therefore $ G_{\m}(M) = G_{\n}(M) $.
 The notation $ G(M) $ will be used
to denote this without any reference to the ring.    Also note
that  $ \depth_{G_{\m}(A)} G(M) = \depth_{G_{\n}(B)}G(M)$.
\end{remark}

 \s Recall that the function $n \mapsto \lambda(M/\m^{n+1}M)$ is called the \emph{Hilbert-Samuel} function.
Let $p_M(z)$ be the \emph{Hilbert-Samuel
polynomial}. The following number
\begin{equation}
\label{postno}
\post(M) = \min\{ n \mid p_M(i) = \lambda(M/\m^{i+1}M) \ \text{for all} \ i \geq n \}
\end{equation}
 is called the \emph{postulation number} of $M$ (\wrt \ $\m$).

\s \label{formulaone} If $f(z) = \sum_{k \geq 0}a_kz^k \in \Z[z]$ then for $i \geq 0$ set $e_i(f) = f^{(i)}(1)/i! = \sum_{k \geq i}\binom{k}{i}a_k$ and set 
\[
\chi_{i}(f) = \sum_{j=0}^{i}(-1)^{i-j}e_{i-j}(f)    + (-1)^{i+1}f(0) = \sum_{k\geq i+1}\binom{k-1}{i}a_k .
\]
It follows that if $a_i \geq 0$ for all $i \geq 0$ then $e_i(f) \geq 0$ and $ \chi_{i}(f) \geq 0$ 
 for all $i \geq 0$. The following Lemma can be easily proved.

\begin{lemma}
\label{formulaonelemma}
 If $g(z), p(z),q(z)$ and $r(z)$ are polynomials with integer coefficients that satisfy
the equation $(1-z)g(z) = p(z) - q(z) + r(z)$ then
\begin{enumerate}[\rm (i)]
 \item
 $e_0(q) = e_0(p) + e_0(r)$. 
\item
 $e_{i}(q) = e_i(p) + e_i(r) + e_{i-1}(g)$ for $i \geq 1$.
 \item 
$\chi_{0}(q) = \chi_{0}(p) + \chi_{0}(r) + g(0)$.
\item
$\chi_{i}(q) = \chi_{i}(p) + \chi_{i}(r) + \chi_{i-1}(g)$ for $i \geq 1$.
\item
If all the coefficients of $g$ are non-negative then for $i \geq 0$ 
we have \\
$e_{i}(q) \geq  e_i(p) + e_i(r)$ and 
$\chi_{i}(q) \geq \chi_{i}(p) + \chi_{i}(r)$. \qed
\end{enumerate}
\end{lemma}

\section{Basic Construction}
\begin{remark}
\label{modstruc}
 For each $n \geq 0$   and $t \geq 0$
 set $L_t(M)_n =   \Tor^{A}_{t}(A/\m^{n+1}, M)$. For $t \geq 0$ let  $L_t(M)
= \bigoplus_{n\geq 0}L_t(M)_n $.
 If $x_1,\ldots,x_s$  is a sequence of elements in  $\m$, then we give
$L_t(M)$ a  structure of a graded
 $A[X_1,\ldots,X_s]$-module  as follows:

 For $ i = 1,\ldots,s$ let
$\xi_i: A/\m^n \rightarrow A/\m^{n+1}$  be the maps given by
   $\xi_i(a + \m^n) = x_ia + \m^{n+1}$. These homomorphisms induces
homomorphisms
 $$ \Tor^{A}_{t}(\xi_i,M) :  \Tor^{A}_{t}(A/\m^n,M) \xar
\Tor^{A}_{t}(A/\m^{n+1},M)$$

 Thus,  for $i = 1,\ldots,s$  and  each $t $ we obtain homogeneous maps of
degree 1:
 $$ X_i : L_t(M) \xar L_t(M) .$$
 For $i,j = 1,\ldots,s$  the equalities $ \xi_i\xi_j =
\xi_j\xi_i$ yields equalities
 $X_iX_j = X_jX_i$.  So $L_t(M)$ is a graded
$A[X_1,\ldots,X_s]$-module for each $t \geq 0$.
\end{remark}

 \begin{proposition}
 \label{propE}
  Let $M, F$, and $K$ be finite $A$-modules and   let $x_1,\ldots,x_s$  be a
sequence of elements in  $\m$.
    If    $L_t(M)$, $L_t(F)$ and $L_t(K)$ are given the
$A[X_1,\ldots,X_s]$-module structure described in Remark
\ref{modstruc}
     then
   \begin{enumerate}[\rm \quad 1.]
   \item
   Every exact sequence of $A$-modules
   $\displaystyle{0 \rt K \rt F \rt M \rt 0 }$
   induces  a long exact sequence of graded $A[X_1,\ldots,X_s]$-modules

   $\displaystyle{ \cdots \rt L_{t+1}(M) \rt L_t(K) \rt L_t(F) \rt L_t(M)
\rt \cdots \rt L_0(M) \rightarrow 0}$.
\item
 For  $i = 1,\ldots,s$    there is an equality
  $$\ker\big(L_0(M)_{n-1} \xrightarrow{X_i}  L_0(M)_{n}\big)  =  \frac{\m^{n+1}M\colon_M x_i}{\m^n M} $$
    \item
    If $x_i \in \m\setminus \m^2$ is such that $x_i^*$ is
    $G(M)$-regular then $X_i$ is $L_0(M)$-regular.
    \item
    If $F$ is free $A$-module and $x_i$ is $K$-superficial for some $i$ then
     \begin{enumerate}
     \item
         $\displaystyle{\ker\big(L_1(M) \xrightarrow{X_i}  L_1(M)\big)_n = 0
\quad \text{for} \ n \gg 0 }$
     \item
         If $x_{i}^{*}$ is $G(K)$-regular then $X_i$ is $L_1(M)$-regular.
    \end{enumerate}
\end{enumerate}
 \end{proposition}

  \begin{proof}
 To prove part 1, set $S = A[X_1,\ldots,X_s]$,
 $$\beta_{t,n} = \Tor^{A}_{t}(A/\m^{n+1},\beta) : \Tor^{A}_{t}(A/\m^{n+1},K)  \xar \Tor^{A}_{t}(A/\m^{n+1},F)  $$
   $$\alpha_{t,n} = \Tor^{A}_{t}(A/\m^{n+1},\alpha) : \Tor^{A}_{t}(A/\m^{n+1},F)  \xar \Tor^{A}_{t}(A/\m^{n+1},M)  $$
    and consider the connecting homomorphisms
    $$\delta_{t+1,n}:  \Tor^{A}_{t+1}(A/\m^{n+1},M) \xar  \Tor^{A}_{t}(A/\m^{n+1},K).$$
 By a well known theorem in Homological algebra
    if $\mathbf{X}$ is a free resolution of $K$ and $\mathbf{Z}$ is a free resolution of $M$
    then there
     exists a  free resolution $\mathbf{Y}$ of $F$ and an exact sequence of
    complexes of free $A$-modules $0 \rt \mathbf{X} \rt \mathbf{Y} \rt \mathbf{Z} \rt 0$
whose homology sequence is the given exact sequence $0 \rt K \rt F \rt M \rt0$.
    This yields for each $n$ a commuting diagram of complexes with exact rows ;
       \[
  \xymatrix
{
 0
 \ar@{->}[r]
  & \mathbf{X}\otimes A/{\m^n}
    \ar@{->}[d]^{\xi_i}
\ar@{->}[r]
 & \mathbf{Y}\otimes A/{\m^n}
    \ar@{->}[d]^{\xi_i}
\ar@{->}[r]
& \mathbf{Z}\otimes A/{\m^n}
    \ar@{->}[d]^{\xi_i}
\ar@{->}[r]
 &0
 \\
 0
 \ar@{->}[r]
  & \mathbf{X}\otimes A/{\m^{n+1}}
\ar@{->}[r]
 & \mathbf{Y}\otimes A/{\m^{n+1}}
\ar@{->}[r]
& \mathbf{Z}\otimes A/{\m^{n+1}}
\ar@{->}[r]
&0
 }
\]

 In homology it induces the following commutative diagram :
    \[
  \xymatrix
{
\cdots
 \ar@{->}[r]
&L_{t+1}(M)_{n-1}
\ar@{->}[d]^{X_i}
\ar@{->}[r]^{\delta_{t+1,n-1}}
&L_{t}(K)_{n-1}
\ar@{->}[d]^{X_i}
\ar@{->}[r] ^{\beta_{t,n-1}}
&L_t(F)_{n-1}
\ar@{->}[d]^{X_i}
 \ar@{->}[r] ^{\alpha_{t,n-1}}
 &\cdots
\\
 \cdots
 \ar@{->}[r]
&L_{t+1}(M)_{n}
\ar@{->}[r]^{\delta_{t+1,n}}
& L_{t}(K)_{n}
\ar@{->}[r]^{\beta_{t,n}}
&L_t(F)_{n}
\ar@{->}[r] ^{\alpha_{t,n} }
 &\cdots
}
\]

This proves the desired assertion.
\begin{remark}
 \label{CommDia}
We will use the exact diagram above often. So when there is a reference to  this remark,
I  mean to refer the commuting diagram above.
\end{remark}
    The second part is clear from the definition of the action $X_i$. Part 3.\ follows from 2.  
    If $F$ is free, then $L_1(F) = 0$, so 1.\
  gives an   exact
    sequence of $S$-modules
$0 \rt L_1(M) \rt L_0(K)$.
Together with 2.\ and 3.\ this
yields the assertions in 4.
 \end{proof}
\begin{remark}
 \label{l1growth}
 If $(A,\m)$ is \CM \ of dimension $d$  and $M$ is maximal non-free \CM \ then by
\cite[Remark 23]{Pu1} there is an equality
 \begin{align}
 \label{*}
 \sum_{n\geq 0}\lambda\left( \Tor^{A}_{1}(M,A/\m^{n+1}) \right) z^n &= \frac{l_M(z)}{(1-z)^d}  \quad
\text{here}\  l_M(z) \in \mathbb{Z}[z] \ \text{and} \ l_M(1) \neq 0 \\
 \label{Beqn}
  (1-z)l_M(z) &= h_{\Syz^{A}_{1}(M)}(z) - \mu(M)h_A(z) + h_M(z).
\end{align}
 \end{remark}

We study the case when $\dim A = 0$.
\begin{lemma}
\label{dim0form}
If $\dim A = 0$ and $M$ is any finite $A$-module then for all $i \geq 0$
\begin{enumerate}[\rm 1.]
\item
$\mu(M)e_i(A) \geq e_i(M) + e_i(\Syz^{A}_{1}(M))$ and $\mu(M)\chi_i(A) \geq \chi_i(M) + \chi_i(\Syz^{A}_{1}(M))$.
\item
$\mu(M)e_i(A) \geq e_i(M)$ and $\mu(M)\chi_i(A) \geq \chi_i(M)$.
\end{enumerate}
\end{lemma}
\begin{proof}
Note that when $\dim A = 0$ we get that $l_M(z)$ has non-negative coefficients. Using
(\ref{Beqn} )  
and Lemma \ref{formulaonelemma}.v  we get 1. and 2. The assertion 3.
follows from 1. and 2.
since $N = \Syz^{A}_{1}(M)$ has dimension zero and so  
$h_{N}(z)$ has non-negative coefficients. Therefore $e_i(N)$ and $\chi_{i}(N)$
are non-negative for $i \geq 0$ (see (\ref{formulaone}) ).
\end{proof}

\s
\label{geqgeq}
It follows from Lemma \ref{dim0form}.3 that 
if $G(A)$ and $G(M)$ is \CM \ then $e_i(A)\mu(M) \geq e_i(M)$ for
all $i \geq 0$.

\begin{remark}
In view of the Remark \ref{l1growth} and (\ref{geqgeq}) it is quite important to understand $L^1(M)= \bigoplus_{n \geq 0}\Tor^{A}_{1}(M,A/\m^{n+1})$ when $M$ is \MCM. In the next Lemma we answer the question when $\dim M = 1$. 
\end{remark}

\begin{lemma}
\label{dim1}
 Let $(A,\m)$ be a Cohen-Macaulay local ring of dimension one,  let $M$
be a non-free maximal Cohen-Macaulay    $A$-modules
 and let
 $$0  \xar E  \xar F \xar M  \xar 0 $$
 be an exact sequence with $F$ a finite free $A$-module.
 Let x be $A\oplus M\oplus E$-superficial. If  $L_1(M)$ is  given
 the
 $A[X]$-module structure  described in Remark \ref{modstruc}
 then
  we have
  \begin{enumerate}[\quad \rm 1.]
  \item
There is an $\m$-primary ideal $\q$ such that $\q A[X]L_1(M) = 0.$ Furthermore
  $L_1(M)$ is a Noetherian
$A/\q[X]$-module of dimension one.
\item
$(1 - z)l_M(z) = h_E(z) - h_F(z) + h_M(z)$.
\item
If $G(E)$ is \CM \ then $X$ is $L_1(M)$-regular. Furthermore
$$e_{i}(F) \geq e_{i}(M) + e_{i}(E)  \quad\text{and} \quad \chi_{i}(F) \geq \chi_{i}(M) + \chi_{i}(E) \ \text{for all $i \geq 0$}.$$ 
\end{enumerate}
\end{lemma}
\begin{proof}
1.  Since $\dim A = 1$  and $M$ is non-free, it follows    from
 Lemma \ref{dim1}
that $\lambda(L_1(M)_n)$ is a non-zero constant for large $ n$.
Since $X \colon L_1(M)_n \rt L_{1}(M)_{n+1}$ is injective for
large $n$ and since $\lambda(L_1(M)_n)$ is constant for large $n$,
it follows that $L_{1}(M)_{n+1} = XL_1(M)_{n}$ for large $n$, say
for all $ n \geq s$.
For $n \geq 0$ set $\q_n = \ann_A L_1(M)_n$.  Note that $\q_n$ is $\m$-primary for all $n$. 
Since the  map
$X \colon  L_1(M)_n \rt L_{1}(M)_{n+1}$ is bijective for all  $n \geq s$ we have  $\q_n = \q_s$ for each $n\geq s$. Set $\q = \cap_{n=0}^{s} \q_n $. Clearly $\q L_1(M)_n = 0$ for each $n\geq 0$. Thus
$L_1(M)$ is an $A/\q[X]$ module.  
 For each $i = 0,1,\ldots, s $ choose a finite
set $\mathcal{P}_i$ of generators  of $L_1(M)_i$ as an $A$-module.
 It is easy to see that $\bigcup_{i=0}^{s} \mathcal{P}_i$
generates $L_1(M)$ over $A[X]$. Since $\lambda( L_1(M)/XL_1(M)) <
\infty$  and   $\lambda( L_1(M)) =    \infty$ it follows that
 $\dim L_1(M) = 1$.

2.\ By Schanuel's lemma, \cite[p.\ 158]{Ma} we have $F\oplus
\Syz^{A}_{1}(M) \cong E\oplus A^{\mu(M)}$. Therefore
$$  (1-z)l_M(z) = h_{\Syz^{A}_{1}(M)}(z) - \mu(M)h_A(z) + h_M(z) =  h_E(z) - h_F(z) + h_M(z). $$

3. It  is clear from 1. and  Proposition \ref{propE}.3 that $X$ is $L_1(M)$-regular.
  It also follows from 1.\ that $l_M(z)$ is the $h$-polynomial of
$L_1(M)$ considered as an $A[X]$-module.
  Set $e^{T}_{i}(M) = l_M^{(i)}(l)/i!$. Since $X$ is  $L_1(M)$-regular
and  $\dim L_1(M) = 1$  we have that
$l_M(z)$ is the  Hilbert series of $L_1(M)/XL_1(M)$. Thus all the
coefficients of $l_M(z)$ is non-negative. Using 2. and Lemma \ref{formulaonelemma}.v 
we get the desired inequalities. 
\end{proof}

\section{Monotonicity}

The following remark will be used often.
\begin{remark}
\label{br2}
 Let $f(z) = \sum_{n\geq 0} a_nz^n $ be a formal power series  with
non-negative coefficients.
If  the power series $g(z) =  \sum_{n\geq 0} b_nz^n $ satisfies $g(z) =
f(z)/(1-z)$, then $b_n = \sum_{i = 0}^{n}a_i $, and so
 the sequence $\verb+{+ b_n \verb+}+$ is nondecreasing.
\end{remark}
The next proposition yields an easy criterion for monotonicity.
\begin{proposition}
\label{PpD}
Let $M$ be an $A$-module. Set $k =A/\m$. If $\depth G(M) \geq 1$ then the Hilbert
function of $M$ is non-decreasing.
\end{proposition}
\begin{proof}
Using Remark \ref{Pr0} we may assume that $k$ is infinite. Thus
  there exists $x \in \m\setminus \m^2$, such that $x^*$
 is $G(M)$-regular.  
  It follows  that the Hilbert function of $M$ is non-decreasing.
\end{proof}

We deduce Theorem   \ref{c1thMo} from the following result.
 \begin{theorem}
  \label{thMo}
  Let $(Q,\m)$ be a local ring with $\depth G(Q) \geq 2$ and let $ M $ be
a  $Q$-module.
  If  $\pd_Q M \leq 1$ then the Hilbert function of $M$ is
non-decreasing.
  \end{theorem}
\begin{proof}
 It is sufficient to consider the case when the residue field of
$Q$ is infinite (see Remark \ref{Pr0}).
  Since $\pd_Q M \leq  1$
we have a presentation of $M$
\begin{equation}
\label{pD1}
  0\xar Q^{n} \xar Q^m \xar M \xar 0 \quad \text{with} \ 0 \leq n \leq m.  
\end{equation}
   Let $x,y$ be elements in $\m \setminus \m^2$ such that  $x^*,y^*$ is a
$G(Q)$-regular sequence. Let  $L_0(Q)$, $L_0(M)$ and $L_1(M)$ be
the   $Q[X,Y]$-modules  described in Remark \ref{modstruc}.

   By Proposition \ref{propE}.3 we get $X$ is $L_0(Q)$-regular.
  Set $B = Q/(x)$ and notice
  $$ \frac{L_0(Q)}{XL_0(Q)} = \bigoplus_{n\geq0}\frac{Q}{(x,\m^{n+1})} =
L_0(B). $$
  Since $G(B) =  G(Q)/x^*G(Q)$ we see that $y^*$ is $G(B)$-regular.
  Proposition \ref{propE}.3 shows that $Y$ is
  $L_0(B)$-regular. Thus $X,Y$ is a $L_0(Q)$-regular sequence.

     Using the exact sequence  (\ref{pD1})  and Proposition \ref{propE}.1, we obtain  an exact sequence of
     graded  $Q[X,Y]$
modules
\begin{equation}
\label{exClaim2}
   0 \xar  L_1(M) \xar L_0(Q)^n \xrightarrow{\phi} L_0(Q)^{m} \xar L_0(M) \xar 0.
\end{equation}
Set $K = \image \phi$. Since $X$ is $L_0(Q)$ regular we see that it is both $K$ and $L_1(M)$-regular. So the exact sequence $0 \rt L_1(M) \rt  L_0(Q)^n \rt K \rt 0$ yields the exact sequence
\[
0 \xar \frac{ L_1(M)}{X L_1(M)} \xar \frac{L_0(Q)^n}{X L_0(Q)^n} \xar \frac{K}{XK} \xar 0
\]
 Since $Y$ is $L_0(Q)^n/X L_0(Q)^n$ regular it follows that $Y$ is $L_1(M)/X L_1(M)$-regular.
Thus  $X,Y$ is an $L_1(M)$- regular sequence.

The regularity of $X,Y$ implies equalities

   \begin{alignat*}{2}
       \sum_{i \geq 0}\lambda(L_0(Q)_i )z^i &= \frac{u(z)}{(1-z)^2}      & \quad \text{where}\ \quad   u(z)  &= \sum_{i \geq 0}\lambda \left(\frac{L_0(Q)_i}{(X,Y)L_0(Q)_{i-1}} \right)z^i.  \\
     \sum_{i \geq 0}\lambda(L_1(M)_i )z^i &= \frac{v(z)}{(1-z)^2}     & \quad \text{where}\ \quad v(z)  &=   \sum_{i \geq 0}\lambda \left(\frac{L_1(M)_i}{ (X,Y)L_1(M)_{i-1}}
     \right)z^i.
   \end{alignat*}

  Using the exact sequence (\ref{exClaim2})
 we get
 \[
\sum_{i\geq 0}\lambda(L_0(M)_i )z^i = (m-n)\frac{u(z)}{(1-z)^2} +
\frac{v(z)}{(1-z)^2}
 \]
The equality $H_M(z) = (1-z)\sum_{i\geq 0}\lambda(L_0(M)_i )z^i $
yields

 $$H_M(z) =  (m-n)u(z)/(1-z) + v(z)/(1-z).$$
  Now Remark \ref{br2} shows   that the Hilbert function of
 $M$ is non-decreasing.
\end{proof}

We obtain Theorem \ref{c1thMo} as a corollary to the  previous theorem.
\begin{proof}[Proof of Theorem  \ref{c1thMo} ]
 We may assume that A is complete and so $A \cong Q/(f)$ for some
regular local ring $(Q, \n)$ and $f \in \n^2$. Then $\depth M = \dim Q - 1$ and
$\pd_Q M = 1$.  Using Theorem \ref{thMo} it follows
 that the Hilbert function of $M$ is non-decreasing.
\end{proof}

Since the Hilbert function is increasing if $\depth G(M) > 0$,
we  construct  a \MCM \ module $M$ over a
hypersurface ring $A$
 such that $\depth G(M) = 0$.

\begin{example}
 Set $Q = k[[ x,y ]] $ and  $\n =(x,y)$. Define $M$ by the
exact sequence
\begin{align*}
0 \xar Q^2 &\xrightarrow{\phi} Q^2 \xar M \xar 0 \quad \text{where}\\
\phi &= \begin{pmatrix}x&y \\-y^2&0\end{pmatrix}
\end{align*}
Set $(A,\m) = (Q/(y^3),\n/(y^3)) $. Note $y^3 = \det(\phi)$
 annihilates $M$. So
   $M$ is a \MCM \ $A$-module.
Set $K = \Syz^{A}_{1}(M)$. Note that $G(Q) = k[x^*,y^*]$.
  Since $y^3M = 0$ ,  we have that if $P \in \Ass_{G(Q)}(G(M))$
  then $P \supseteq  (y^*)$. So we get that $x^* \notin P$ if $P$
  is a relevant associated prime  of $G(M)$. Therefore $x^*$ is an
   $M$-superficial element.
We  show    $ \depth G(M) = 0$. 
 Otherwise by \ref{Pr2}.5.a  we get 
  that $x^*$ is  $G(M)$-regular. However if $m_1,m_2$
are the generators of $M$ then $xm_1 = y^2m_2 \in \n^2M$ and this
 implies  $m_1 \in (\n^2M\colon_M x) = \n M$, which is a
contradiction.
\end{example}

The next corollary
partly overlaps with a result of Elias \cite{ElMo}: all
   equicharacteristic \CM \ rings of dimension 1 and embedding
   dimension 3 have non-decreasing Hilbert functions.
See \cite[p.\ 337]{VaSix} for  an example of a 
 complete intersection ring
    $(A,\m)$  of   dimension 1 and
codimension $2$ such that  $\depth G(A) = 0$.
\begin{corollary}
   \label{c2thMo}
    If $(A,\m)$ be a  complete intersection of  positive dimension and
codimension $2$ then the Hilbert function of $A$ is
    non-decreasing.  
    \end{corollary}
\begin{proof}
We may assume that $A$ is complete and hence $A = Q/(f,g)$ for a
regular sequence $f,g$ in a regular local ring $(Q, \q)$. Set
$(R,\n) = (Q/(f), \q/(f))$.
 Then $G(R)$ is Cohen-Macaulay, $\dim A = \dim R - 1$ and $\pd_R A = 1$.
 Therefore by  Theorem \ref{thMo} we get that the Hilbert function of $A$ is
non-decreasing.
\end{proof}

Another application of  Theorem \ref{thMo} yields the following:
\begin{theorem}
\label{muI}
Let $(A,\m)$ be a Noetherian equicharacteristic local ring of dimension $d > 0$ and let 
$M$ be a \MCM \ $A$-module. Let $I$ be an $\m$-primary ideal in $A$ with $\mu(I) = d +1$.
Then the Hilbert function of $M$ \wrt \ $I$ is non-decreasing.
\end{theorem}
\begin{proof}
Without any loss of generality we may assume that $A$ is complete.  
Let $I = (x_1,\ldots,x_{d+1})$. Since $A$ is complete and equicharacteristic it contains
 a subfield $k \cong A/\m$. Set $R = k[[T_1,\ldots,T_{d+1}]]$
 and let $\n$ be its unique maximal ideal.
 Consider the local homomorphism $\phi : R \rt A$ defined by $\phi(T_i) = x_i$. Then
$A$ becomes an $R$-module via $\phi$. Since $A/\n A = A/I$ has finite length we get 
 $M$ is  a finite $R$-module. It can be easily checked that
 $M$ is  a \CM \ $R$-module of 
dimension $d$.

Since $R$ is regular, $\projdim_R M $ is finite. So
$\projdim_R M = \depth R - \depth M = 1$. Therefore by Theorem \ref{thMo} it follows that
 $H_{\n}(M,j)$ is non-decreasing. 
Note that $\n^j M = I^jM$ for each $j \geq 1$ and so   $H_{\n}(M,j) = \lambda(I^nM/I^{n+1}M)$ for
each $j \geq 0$. 
This establishes the assertion of the theorem.
\end{proof}

\section{ Hilbert coefficients}

\s In this section $\varepsilon_s$ denotes the $s\times s$ identity matrix.
Let $(Q,\n)$ be a regular local ring, $f \in \n^e \setminus\n^{e+1}$,
 $e \geq 2$, $A= Q/(f)$,  $M$ a \MCM \ $A$-module
and $K = \Syz^{A}_{1}(M)$.

 By a \emph{matrix-factorization} of $f$
we mean a pair $(\phi,\psi)$ of square-matrices with elements in $Q$ such that
$$   \phi \psi = \psi \phi = f\ve. $$
If $M$ is an $A$-module then $\pd_Q M = 1$. Also a presentation of $M$
\[
0 \xar Q^n \xrightarrow{\phi} Q^n \xar M \xar 0
\]
 yields a matrix factorization of  $f$.
See \cite[p.\ 53]{EisMf} for details.

In the sequel $(\phi_M,\psi_M)$ will denote a matrix factorization of $f$ such that
\[
0 \xar Q^n \xrightarrow{\phi_ M}  Q^n \xar M \xar 0
\]
is a minimal presentation of $M$.
 Note that
 \[
0 \xar Q^n \xrightarrow{\psi_ M}  Q^n \xar \Syz^{A}_{1}(M) \xar 0
\]
is a not-necessarily minimal presentation of $\Syz^{A}_{1}(M)$.

If $\phi : Q^n \xar Q^m $ is  a linear map then
we set
\begin{equation*}
 i_{\phi} = \max\{i  \mid \ \text{all entries of $\phi$ are in } \  \n^i \}.
\end{equation*}
If
$M$  has  minimal presentations:
$ 0 \rt Q^{n} \xrightarrow{\phi}Q^{n} \rt M \rt 0$ and

\noindent $ 0 \rt
Q^{n} \xrightarrow{\phi'}Q^{n} \rt M \rt 0$,
   then  it is well known that  $i_{\phi} = i_{\phi'}$ and
 $\det(\phi) = u\det(\phi')$ with $u$ a unit.
We set $i(M)= i_{\phi}$ and $\det(M) = \left(\det(\phi)\right)$.
For $g \in Q$, $g \neq 0$,  set $v_Q(g) = \max \{i \mid g \in \n^{e} \}$. For convenience
set $v_Q(0) = \infty$.
Note that $e(Q/(g)) = v_Q(g)$ for any $g \neq 0$.
We first consider the  case when $\dim A = 0$.

\begin{remark}
\label{dim0}
Let $(Q,\n)$ be a DVR, $v_Q(f) = e$, $A = Q/(f)$ and $M$ a finite $A$-module. If
$\n = (y)$ then $f = uy^e$,
where $u$ is a unit.
   Therefore as an $Q$-module
   $$M \cong \bigoplus_{i=1}^{\mu(M)}Q/(y^{a_i}) \quad \text{with} \ 1 \leq a_1 \leq \ldots \leq a_{\mu(M)} \leq e.$$

   This  yields a minimal presentation of $M$:
  $$ 0 \rt Q^n \xrightarrow{\psi} Q^n \rt M \rt 0  \ \text{where}\ \psi_{ij} = \delta_{ij}y^{a_i}.            $$
   This yields
\begin{enumerate}
\item
$i(M) = a_1$.
\item
$h_M(z) = \mu(M)(1+z+\ldots + z^{i(M) - 1}) + $ higher powers of $z$
\item
$h_0(M)\geq h_1(M) \geq \ldots \geq h_s(M)$.
\item
$M$ is free if and only if $i(M) = e$.
\item
As an $Q$-module
$K \cong \bigoplus_{i=1}^{\mu(M)}Q/(y^{e - a_i})$.
\item
     $e(M) \geq \mu(M)i(M)$ and $e_1(M) \geq \mu(M)\binom{i(M)}{2}$.
\item
$v_Q(\det \phi) \geq i(M)\mu(M)$ with equality iff $e(M) = i(M)\mu(M)$.
\end{enumerate}
\end{remark}

We note an immediate corollary to  assertion 3. in the  previous remark.
\begin{corollary}
Let $(A,\m)$ be a hypersurface ring of dimension $d$. Let $M$ be a \MCM \ $A$-module such that $G(M)$ is \CM.
If $h_M(z) = h_0(M) + h_1(M)z + \cdots + h_s(M)z^s$  is the $h$-polynomial of $M$
 then $h_0(M)\geq h_1(M) \geq \ldots \geq h_s(M)$.
\end{corollary}
\begin{proof}
 We may assume that $A$ is complete with
infinite residue field, hence $A = Q/(f)$ for some regular local ring
$(Q,\n)$. Consider $M$ as a $Q$-module. Let $x_1,\ldots,x_d$ be a
$Q\oplus M$-superficial sequence. Set $J = (x_1,\ldots,x_d)$,
$(R,\q) = (Q/J,\n/J)$, $B = A/J$ and $N = M/JM$. Note that $R$ is
a DVR. Since $G(M)$
is \CM \ we also have $h_M(z) = h_{N}(z)$ and so the result
follows from Remark \ref{dim0}(3).
\end{proof}

To use the other assertions in Remark \ref{dim0} we need the following
definitions. The notion of superficial sequence is extremely  useful in the study of
Hilbert functions. We need to generalize it to deal also with
 a presentation of 
a module.

\begin{definition}
Let $(Q,\n)$ be a regular local ring, $f \in \n^e \setminus\n^{e+1}$, $e \geq 2$, $A= Q/(f)$ and
  $M$ a \MCM \ $A$-module.
Let $ 0 \rt Q^{n} \xrightarrow{\phi}Q^{n} \rt M \rt 0$ be a minimal
presentation of $M$.  We say that $x \in \n$ is $\phi$-\emph{superficial} if
\begin{enumerate}
\item
$x$ is $(Q\oplus M \oplus A)$-superficial.
\item
If $\phi = (\phi_{ij})$ then $v_Q(\phi_{ij}) = v_{Q/xQ}(\ov{\phi_{ij}})$
\item
$v_Q(\det(\phi)) = v_{Q/xQ}(\det(\ov{\phi}))$.
\end{enumerate}
\end{definition}
Since $e(Q/(g)) = v_Q(g)$ for any $g \neq 0$ it follows that
 if  $x$ is $Q\oplus M \oplus A \oplus \left(\bigoplus_{ij}Q/(\phi_{ij})\right) \oplus Q/\det(\phi )$ -superficial
then it is $\phi$-superficial. So
$\phi$-superficial elements exist
 if the residue field of $Q$ is infinite.

 If $x$ is $\phi$-superficial, then clearly
$i(M) = i(M/xM)$. Also note that
$Q/xQ$ is regular and we have an exact sequence
\[
 0 \xar \left(\frac{Q}{xQ}\right)^{n} \xrightarrow{\phi\otimes_Q Q/xQ}\left(\frac{Q}{xQ}\right)^{n} \xar \frac{M}{xM} \xar 0
\]
 This enables the following definition:
\begin{definition}
Let $(Q,\n)$ be a regular local ring, $f \in \n^e \setminus\n^{e+1}$, $e \geq 2$, $A= Q/(f)$ and
  $M$ a \MCM \ $A$-module. Let $ 0 \rt Q^{n} \xrightarrow{\phi}Q^{n} \rt M \rt 0$ be a minimal
presentation of $M$.  We say that $x_1,\ldots,x_r$ is a
 $\phi$-\emph{superficial
sequence}  if
$\ov{x_i}$ is $\left(\phi\otimes_Q Q/(x_1,\ldots,x_{i-1})\right)$-superficial for $i = 1,\ldots,r$.
\end{definition}

\noindent \textbf{Notation:} Let $M$ be  an $A$- module. If $x$  is  $A\oplus M$ superficial (or more generally
 it is  superficial  \wrt \ to an injective  map $\theta :Q^n \rt Q^n $)
 then set  $(B,\n) = (A/(x),\m/(x)) $ and $N = M/xM $.

We need a few preliminaries  before
 we prove Theorem 2.
\begin{lemma}
\label{U1}
 Let $(A,\m)$ be a \CM \ local ring of dimension $d >0$ with
 infinite residue field.
Let $M$ be a \CM \ $A$-module of dimension $1$ with a presentation
$G \xrightarrow{\phi}F \xar M \xar 0 $  such that all entries in
$\phi$ are in $\m^{l}$. If $\depth G(A) \geq 1$ and $x$ is a
$A\oplus M$-superficial element  then
\begin{enumerate}[\rm \quad 1.]
\item
$(\m^{i+1}M\colon_M x) = \m^iM  $ for $i = 0,\ldots,l-1 $.
\item
Furthermore if   $\m^{l}M \subseteq xM$ then $\depth G(M) \geq 1$.
\end{enumerate}
\end{lemma}
\begin{proof}
Set $b_i(M) = \lambda((\m^{i+1}M\colon_M x)/\m^iM) $. Since all
the entries of $\phi$ are in $\m^l$ we have that $\phi_{0,j-1} =
\phi\otimes A/\m^{j} = 0$ for $j = 1,\ldots,l$.

1. Note that $b_0(K) = 0$ for any $A$-module $K$.  Fix $ i $ with $
1\leq i \leq l-1$.
 Let $p \in (\m^{i+1}M \colon_M x)$.
 Let $u $ be the pre-image of $p$ in $F$. Using the
commutative diagram  \ref{CommDia} and since $ \phi_{0,i} = 0$ and $ \phi_{0,i+1} = 0$ we
obtain that $xu \in \m^{i+1}F$. Since $x$ is $A$-superficial and
depth $G(A) \geq 1$ we get by  \ref{Pr2}.5.a that $x^*$ is
also $G(A)$-regular. Therefore $ u \in \m^iF$ and so $p \in
\m^iM$. Thus $b_i(M) = 0$.

2.  Since
$\m^lM \subseteq xM$ we get that $\n^lN = 0$ and so
$\sum_{i=0}^{l-1} H(N,i) = e(N) = e(M)$.  Using  \ref{Pr2}.4.b we get that
\[
H(M,l-1) = \sum_{i=0}^{l-1} H(N,i) - b_{l-1}(M) = e(M)
\]
Using  \ref{Pr2}.6.a we obtain $\m^{l}M =x\m^{l-1}M$. So we
have that $\m^{i +1}M =x\m^{i}M$ for all $i \geq l-1$. So we
obtain that $b_i(M) = 0$ for all $i \geq l-1$. This combined with
1.  yields that $x^*$ is $G(M)$-regular.
\end{proof}
An interesting consequence of the  lemma above is the following
lemma  which gives information about the Hilbert function
of a \MCM \ module over a hypersurface ring of dimension 1.
\begin{lemma}
\label{dim1a}
Let $(Q,\n)$ be a regular local ring of dimension two, 
 $f \in \n^e \setminus\n^{e+1}$, $e \geq 2$, $A= Q/(f)$.
If  $M$ is  a \MCM \ $A$-module,
then
\[
h_M(z) = \mu(M)(1 + z + \ldots + z^{i(M) - 1}) + \sum_{i \geq i(M)}h_{i}(M)z^i \quad \text{and} \ h_i(M) \geq 0 \ \forall   i.
\]
\end{lemma}
\begin{proof}
As $\dim M = 1$, the Hilbert series of $M$ 
is  $h_M(z)/(1-z)$.
The Hilbert function of $M$ is non-decreasing by Theorem \ref{c1thMo}.
Therefore all the coefficients of $h_M(z)$ are non-negative.
Set $b_i(M) = \lambda((\m^{i+1}M\colon_M x)/\m^iM) $.
Since
$$A^n \xrightarrow{\phi \otimes A} A^n  \xar M \xar 0 $$
is exact and all the entries of $\phi$  are in $i(M)$ we get
by Lemma \ref{U1}.1 that $b_i(M) = 0 $ for $i = 0,\ldots,i(M)-1$.
This and Remark 4.1.1 yields that
\[
h_M(z) = \mu(M)(1 + z + \ldots + z^{i(M) - 1}) + \sum_{i \geq i(M)}h_{i}(M).
\]
\end{proof}
Next we get an upper bound on $l$ such that $\m^{l}M
\subseteq xM$ holds.
\begin{remark}
\label{Ur1}
 If $\dim A = 1$ and $x$ is $A$-superficial then note that
 since the ring $B$ has length $e_0(A)$ we get that $\m^{e_0(A)} \subseteq (x)$. Therefore if $\dim A = 1$,
 $M$  a maximal $A$-module and  $x$ is $A \oplus M$-superficial  then     $\m^{e_0(A)}M \subseteq (x)M$
\end{remark}

The next lemma deals with the case when $M$ is a syzygy of a
\MCM \ $A$-module.

 \begin{lemma}
\label{U3}
 Let $(A,\m)$ be a \CM \ $A$-module of dimension $1$and let $L$ be a
  non free \MCM \  $A$-module. Set $M =
\Syz^{A}_{1}(L)$. If
$x$  is $(A\oplus M \oplus L)$-superficial
then
 $\m^{e _0(A)-1 }M \subseteq ( x)M$.
\end{lemma}
\begin{proof}
 By Remark
\ref{Ur1} we have $\m^{e _0(A) } \subseteq (x)$. We
also have an exact sequence : $0 \rt M \rt F \rt L \rt 0$
where $F$ is a free $A$-module.
Set  $G = F/xF$ and $W = L/xL$.
Going mod $x$ we get $0 \rt
N \rt G \rt  W \rt 0$. Note that $N
\subseteq \n G$. Therefore  $\n^{e _0(A)-1 }N
\subseteq \n^{e _0(A) }G = 0$. It follows that
$\m^{e _0(A)-1 }M \subseteq xM$.
\end{proof}

\begin{proof}[Proof of Theorem \ref{mtHy}]
Clearly we may assume that  $k = Q/\n$ is infinite.
 Let  $ 0 \rt Q^n \xrightarrow{\phi} Q^n  \rt M \rt 0$ be
 a minimal presentation of $M$ over $Q$.
 If $\dim A \geq 2$, then choose
$x_1,\ldots,x_d$ to be a  maximal $\phi$-superficial sequence.
 Set $J = (x_1,\ldots,x_{d-1})$. Since all the invariants considered in the theorem remain same
 modulo $J$ it suffices to assume  $\dim A \leq 1$.
When $\dim A = 0$ then all  the results follow easily
by Remark \ref{dim0}.

Therefore assume that
$\dim A = 1$.
Let $x$ be $\phi$-superficial.
Set
 $R = Q/xQ$, $N = M/xM$, $\ov{f} = $ the image of $f$ in $R$ and $B = A/xA = R/(\ov{f})$.
  Note that
  \begin{enumerate} [\rm a.]
  \item
   $i(M) = i(N)$.
   \item
 $e(M) = e(N)$ and $e_1(M) \geq e_1(N)$   ( by 1.2.3 and 1.2.4 )
 \item
 $  R$ is a DVR with maximal ideal say $\q =  (y)$.
  \item
  $v_{R}(f) = v_{Q}(f) = e $.
   \end{enumerate}

   So 1. follows from Remark \ref{dim0}  and (b) above.

   2. If $i(M) = e$  then  since $i(M) = i(N)$ we get
by Remark \ref{dim0}.4 that $N$ is a free $B$ module.  So
   $M$ is a free $A$-module.  Conversely 
 if $M$ is free then clearly $i(M)= e$.

3. Let $M = F\oplus L$
  where $F$ is a free $A$-module and $L$ has no free summands. Note that $i(L) = i(M)$.
Since $G(A)$ is \CM \ it suffices to show $G(L)$ is \CM.  Notice $L =\Syz^{A}_{1}\left(\Syz^{A}_{1}(L)\right)$.
 If
$x$  is a $A\oplus L \oplus N$-superficial, 
element, then by Lemma \ref{U3} we get
 $\m^{e _0(A)-1 }L \subseteq xL$.
Since $A^n \rt A^n \rt L \rt 0$ is exact  and $\depth G(A) = 1$ we get by Lemma \ref{U1}.2 that $G(L)$ is \CM.

4.  By Proposition \ref{dim1a} we get that
$$h_M(z) = \mu(M)(1 + z + \ldots + z^{i(M) - 1}) + \sum_{i \geq i(M)}h_{i}(M)z^i \quad \text{and}\ h_i(M) \geq 0 \  \forall i.  $$
 It follows that (i) and (ii) are equivalent.
The assertion (iii) $\implies$ (ii) is clear.

(i) $\implies$ (iii). Note that  $\mu(N) = \mu(M)$ and

$$h_N(z) = \mu(M)(1 + z + \ldots + z^{i(M) - 1}) +
 \sum_{i \geq i(M)}h_{i}(N)z^i.$$

Also all the coefficients are non-negative.
Therefore 
$e(M)= e(N) = i(M)\mu(M)$
 if and only if
$h_M(z) = h_N(z) = \mu(M)(1 + z + \ldots + z^{i(M) - 1})$.
Since $h_M(z) = h_N(z)$ we also get  that $G(M)$ is \CM \ (see \ref{Pr2}.8).

Note that since
$M$ is not free $i(M) < e$.
We first assert that $M$ has no free summands.
Otherwise $M = F\oplus W$ where $F$ is free.
This yields $h_M(z) = h_F(z) + h_W(z)$. Since all the coefficients of
$h_F(z)$ and $h_W(z)$ are non-negative
we get that coefficient of $z^{e-1}$  is non-zero. This contradicts
(c).
Therefore if $(\phi,\psi)$ is a matrix-factorization of $M$
then we have a minimal presentation of
$K$ 
\[
0 \xar Q^n \xrightarrow{\psi} Q^n \xar K \xar 0.
\]
Let $x$ be both $\phi$ and $\psi$-superficial.
Set $N = M/xM$ and $R = Q/(x)$.
Since (a) holds then note that $N \cong \big(R/(y^{i(M)})\big)^{\mu(M)}$.
Then
\[
\Syz^{B}_{1}(N) \cong \big(R/(y^{e -i(M)})\big)^{ \mu(M)}.
\]
Since $\Syz^{B}_{1}(N)  \cong K/xK$ we get that
$i(K) = e -\mu(M)$ and so
 $$e(K) = e(K/xK) = \mu(K/xK)i(K) = \mu(K)i(K).$$
Therefore  
by the equivalence of (i) and (iii) we get the required result.
\end{proof}

 \begin{remark}
 \label{U}
   Theorem  \ref{mtHy} can be applied  to  the case of Ulrich modules, that is,
    \MCM \  modules that satisfy $e(M) =
\mu(M)$.
  It is known, see \cite{HUB},  that  Ulrich
$A$-modules exist  when $A$ is a complete hypersurface ring.  Using the previous theorem we get that if $M$ is Ulrich, then
$i(M)= 1 $ and so $G(\Syz^{A}_{1}(M))$ is \CM. Furthermore $i(\Syz^{A}_{1}(A)) =e -1$ and 
$h_{\Syz^{A}_{1}(M)} = \ \mu(M) (1+z+\ldots+ z^{e-2}).$
 \end{remark}

An easy way to test the hypothesis of the previous theorem in the equicharacteristic case is the following:
\begin{proposition}
\label{eqchar}
Let $Q = k[[y_1,\ldots,y_{d+1}]]$.
Let $M$ be a $Q$-module with a minimal presentation
$0 \rt Q^n \xrightarrow{\phi}Q^n \rt M \rt 0$. Set
\[
\phi = \sum_{i \geq i(M)}\phi_i \quad \text{where $\phi_i$ are forms of degree $i$}
\]
Then  $\det \phi_{i(M)} \neq 0$
\IFF \  $h_M(z) = \mu(M)(1+z+\ldots+ z^{i(M)-1})$.
\end{proposition}
\begin{proof}
Note that $\det \phi_{i(M)} \neq 0$ \IFF \ $v_Q(\det \phi) = i(M)\mu(M)$.

Let $f = \det \phi$. Note that $M$ is a maximal $A = Q/(f)$-module.
Let $x_1,\ldots,x_d$ be a maximal $\phi$-superficial sequence.
Set $J = (x_1,\ldots,x_d)$, $R = Q/J$, $\ov{f} = $ image of
$f$ in $R$, $N = M/JM$ and $\ov{\phi} = \phi \otimes Q/J$. Note
 that
$$ i(N) = i(M)\quad \text{and} \quad  v_R(\det \ov{\phi}) = v_Q(\det \phi).$$
 If $\det \phi_{i(M)} \neq 0$ then $v_Q(\det \phi) = i(M)\mu(M)$.
 So  $v_R(\det \ov{\phi}) = i(M)\mu(M) $.
Therefore by Remark \ref{dim0}.7 we get that $e(N) = \mu(N)i(N)$.
This yields  $e(M) = i(M)\mu(M)$ and so by Theorem \ref{mtHy} we get the required
assertion.
Conversely if  $h_M(z) = \mu(M)(1+z+\ldots+ z^{i(M)-1})$ then by
Theorem \ref{mtHy} we get that $G(M)$ is \CM.
So $h_N(z) = h_M(z)$. So we get  $e(N) = \mu(N)i(N)$.  Therefore by Remark \ref{dim0}.7
 we get $v_R(\det \ov{\phi}) = i(N)\mu(N) = i(M)\mu(M) $. So   $v_Q(\det \phi) = i(M)\mu(M)$.
\end{proof}

   We  give an application of the proposition  proved above.
\begin{example}
\label{f}
Set $Q = k[[y_1,\ldots,y_{d+1}]]$ and $\n$ to be the maximal ideal of $Q$. Let  $a,b,c,d$ be in
$\n$ be such that $ f = ad -bc  \neq 0$. Set $A = Q/(f)$.
Set
\begin{gather*}
\phi = \begin{pmatrix}a & b \\c&d\end{pmatrix}\quad \text{and} \quad \psi = \begin{pmatrix}d&-b \\-c&a \end{pmatrix}
\end{gather*}
Define $M$ and $K$ by the exact sequences
\[
0 \rt Q^2 \xrightarrow{\phi}  Q^2 \rt M \rt 0 \quad \text{and} \quad 0 \rt Q^2 \xrightarrow{\psi} Q^2 \rt K \rt 0
\]
Note that $M$ and $K$ are $A$-modules and $K = \Syz^{A}_{1}(M)$.
\begin{enumerate} [\rm  1.]
\item
If $f \in \m^2 \setminus \m^3$ then $M$ is an Ulrich $A$-module.
\item
If $f \in \m^3 \setminus \m^4$ then $M$ or $K$ has minimal multiplicity. Both
are not Ulrich.
\end{enumerate}
To prove 1. note that $i(M) = 1$ and $\det(\phi_1) \neq 0$. So by Proposition \ref{eqchar}
we get that $h_M(z) = \mu(M)$. So $M$ is Ulrich.

2. We first show that $M$ and $K$ are not Ulrich. As $f \in \n^3 \setminus \n^4$
we  have $i(M) = 1$. Since $f \in \n^3$ we also get $\det(\phi_1) = 0$.
  By Proposition \ref{eqchar}
we get that $h_M(z) \neq \mu(M)$. So $M$ is not Ulrich. Similarly we get $K$
is not Ulrich.
Notice $h_A(z) = 1 +z + z^2$. By (\ref{Beqn}) and Proposition \ref{formulaone}.(iv) we have
$$\mu(M)\chi_1(A) - \chi_1(M) - \chi_1(K) = \chi_0(l_M(z)) = l_M(1) - l_M(0).$$
Notice $l_M(0) = \mu(K)$. 
By \cite[Lemma 19]{Pu1} $l_M(1) \geq e_0(K)$. Since $K$ is not Ulrich we have
$\chi_0(l_M(z)) > 0$. It follows
that  $\chi_1(M) = 0$ or $\chi_1(K) = 0$.  By \ref{Pr2}.9 we
get that $M$ or $K$ has minimal multiplicity.
\end{example}

 \section{Second Hilbert coefficient}

\s\label{hypothesis}
In this section $A$ is \CM \  and $M$ is a \MCM \ $A$-module. Set $k = A/\m$.

In view of \ref{cricket1} is natural to ask how do higher Hilbert coefficients of $A$, $M$ and $\Syz^{A}_{1}(M)$ are related.
Lemma  \ref{dim1} indicates a way.

\begin{theorem}
\label{e2general}
  (with hypothesis as in \ref{hypothesis}) Assume 
we have an exact sequence $0 \rt M \rt F \rt E \rt 0$ with $F$ free $A$-module and
 $E$  a finite 
\MCM \ $A$-module. If $G(M)$ is \CM \ and $\depth G(A) \geq d-1$ then
\begin{enumerate}[\rm 1.]
\item
$\displaystyle{e_2(F) \geq e_2(M) + e_2(E) \quad \text{and} \quad \chi_2(F) \geq \chi_2(M) + \chi_2(E).}$ 
\item
$\displaystyle{e_i(F) \geq e_i(M)  \quad \text{and} \quad \chi_i(F) \geq \chi_i(M)}$ for $i \geq 0$. 
\end{enumerate} 
\end{theorem}
  Theorem \ref{e2general} is not satisfactory as there is no easy criteria for 
finding an \MCM \ $A$-module $E$ with   $G(\Syz^{A}_{1}(E))$ is \CM.
However if $A$ is \emph{Gorenstein} then every \MCM \ $A$-module is a syzygy of a \MCM \ $A$-module.
This can be  seen as follows:
Let $M$ be a \MCM \ $A$-module.
 The module $M^* = \Hom_A(M,A)$ is also a \MCM \ $A$-module. Let $N = \Syz^{A}_{1}(M^*)$ and
 $F = A^{\type(M)}$. We 
have an exact sequence $0 \rt N \rt F \rt M^* \rt 0$. Dualizing we get
$0 \rt M^{**} \rt F^* \rt N^* \rt 0.$
Note that $M \cong M^{**}$. Set $S^A(M) =N^*$. Interestingly  $S^A(M)$ behaves
well  mod superficial sequences.

\begin{lemma}
\label{l1the2}
(with hypothesis as in \ref{hypothesis}). Let $A$ be a Gorenstein ring.
   If $x$ is $A\oplus M$-regular
then for  the $B = A/(x)$-module $N
= M/xM $ we have
\begin{enumerate}[\quad \rm 1.]
\item
$\Hom_{B}(N,B) \cong M^*/xM^*$.
\item
$\Syz^{B}_{1}(M^*/xM^*)  \cong
\Syz^{A}_{1}(M^*)/x\Syz^{A}_{1}(M^*)$.
\item
$S^{B}(N) \cong S^A(M)/xS^A(M)$.
\end{enumerate}
\end{lemma}

\begin{proof}
1. We use 
$ 0 \rt M \xrightarrow{x} M \rt N \rt 0 $
to get a long exact sequence
\begin{align*}
0 &\xar \Hom_A(N,A) \xar \Hom_A(M,A) \xrightarrow{x} \Hom_A(M,A) \\
  &\xar  \Ext_{A}^{1}(N,A) \xar \Ext_{A}^{1}(M,A) .
\end{align*}
Notice  $\Ext_{A}^{1}(M,A)= 0 $ as $M$ is \MCM \ and $A$ is Gorenstein \cite[3.3.10.d]{BH}. Using  the
isomorphisms
$ \Hom_A(N,A) = 0$ and
$\Ext_{A}^{1}(N,A) \cong \Hom_{B}(N,B) $
 (see \cite[3.1.16]{BH}),
we get 1.  Note that 2. holds since $M^*$ is
\MCM. To prove 3. we use   1. to get
\[
S^A(M)/xS^A(M) =  \frac{\Hom_A\big(\Syz^{A}_{1}(M^*),A\big)}{x\Hom_A\big(\Syz^{A}_{1}(M^*),A\big)} 
  \quad \cong\Hom_{B}\left(\frac{\Syz^{A}_{1}(M^*)}{x\Syz^{A}_{1}(M^*)}, B\right)
\]
  \[
  S^{B}(N) =    \Hom_{B}\big( \Syz^{B}_{1}(\Hom_{B}(N,B)),B\big) \quad 
                \cong  \Hom_{B} \big(\Syz^{B}_{1}(M^*/xM^*),B\big) 
\]
Finally we use 2. to get the result.
\end{proof}

The following example shows that the hypothesis on 
$\depth G(A)$  in Theorem \ref{e2general} cannot be dropped.
\begin{example}
 Set
$R = k[[x,y,z,u,v]]$ and
 $\q = (z^2,zu,zv,uv,yz-u^3,xz-v^3)$. Set
$A = R/\q$, $E = R/\q + (z)$ and $M = (\q + (z))/\q$. 
The ring $A$ is \CM \ and by \cite{cocoa} we get
  $h_A(t) = 1+ 3t + 3t^3 - t^4$. It is known $\depth G(A) = 0$, see \cite[3.10]{math.AC/0304111}.
 Note that $M$ and
 $E$ are  $A$-modules and we have
an obvious exact sequence $0 \rt M \rt A \rt E \rt 0$. We show
\begin{enumerate}[\rm (a)]
\item
$E$ and $M$ are \MCM \ $A$-modules and $G(E)$ and $G(M)$  are \CM.
\item
$\mu(E)e_3(A) \ngeq e_3(M) + e_3(E)$ and 
$\mu(E)e_3(A) \ngeq e_3(M)$.
\end{enumerate}

\noindent\textit{Proof} (a). If we prove $E$ is \MCM \ then it follows that $M$ is \MCM.  Notice 
$$E = \frac{k[[ x,y,u,v]]}{ (z,uv,u^3,v^3)} \cong \frac{k[[ x,y,u,v]]}{(uv,u^3,v^3)} \cong \frac{ k[[u,v]]}{(uv,u^3,v^3)}\big{[[x,y]]}. $$ 
So $E$ is \MCM. Clearly
$G_{\m}(E) = \frac{ k[[u,v]]}{(uv,u^3,v^3)}\big{[x^*,y^*]}.$
So $ G(E)$ is \CM \ and $h_E(t)  = 1+ 2t + 2t^2$.
 Notice $M$ is  a cyclic  $A$-module.
 Since
$(z,u,v)M = 0$ we get that $M$ is a $ S = A/(z,u,v) = k[[x,y]]$ module. Since $M$ is also 
\MCM \ $S$-module it is free.  As $M$ is cyclic and free $M \cong S$.
Thus $G(M) \cong k[x^*,y^*]$ is \CM \ and  $h_{M}(t) = 1$.

(b) Note that $e_3(A) = -1$, $\mu(E) = 1$, $e_3(M) = e_3(E) = 0$.
\end{example}
\begin{proof}[Proof of Theorem \ref{e2general}]
We  prove the result regarding $e_i$. The  result regarding 
$\chi_i$ can be proved on similar lines.
By Remark \ref{Pr0} we may assume $k$
is infinite.

     When $\dim A = 1$,  assertion 1. follows from Lemma
     \ref{dim1}.
 When $\dim M = 2$, let $x$ be  $M \oplus E \oplus A$-superficial.
  Set     $(B,\n) = (A/(x),\m/(x))$,       $N = M/xM$ and $G = F/xF$ and $L = E/xL$.
Note that $B$ is one-dimensional \CM \ ring and $N$, $L$ are  \MCM \ $B$ modules,
and we have an exact sequence $0\rt N \rt G \rt L \rt 0$.

  Since $G(M)$ is Cohen-Macaulay and $\depth G(A) \geq 1$
  we have that $e_2(M) = e_2(N)$ and $e_2(F) = e_2(G)$.
 Furthermore it follows from \ref{Pr0}.6 that
 $e_2(E) \leq e_2(L)$.
 Therefore we have:
\[
 e_2(M) + e_2(E) \leq e_2(N) + e_2(L)  \leq   e_2(G) = e_2(F).
\]
 Note that the second inequality above follows from the dimension one case.

   When $\dim A > 2$ let $x_1,\ldots,x_{d-2}$ be a $M \oplus E
   \oplus A$-superficial sequence.
    Set $J = (x_1,\ldots , x_{d-2}), (B,\n) = (A/J, \m/J), N = M/JM, G = F/JF $ and $L = E/JE$.
    By \ref{Pr0}.6 we get that $e_2(E) =
e_2(L), e_2(M) = e_2(N)$
     and $e_2(F) = e_2(G)$.
     So the result follows from the dimension 2 case.

       To prove (ii) note that since $G(M)$ is Cohen-Macaulay and $\depth
G(A) \geq d - 1$
       it suffices to consider the case when $\dim A = 1$.
       By \ref{dim1}.4   we get that
$e_i(S^A(M)) \geq 0$ for all $i \geq 1$.
        So we get $e_{i}(A)\mu(M) \geq e_{i}(M)$ by Lemma
\ref{dim1}.
\end{proof}

Theorem  \ref{the2} now follows as a corollary to  Theorem  \ref{e2general}.

    \begin{proof}[Proof of Theorem \ref{the2}]
By Remark \ref{Pr0} we may assume $k$
is infinite.   Set   $L^* = S^{A}(M)$ and $ F^* = A^{\type(M)}$.
We have an exact sequence
 $0  \rt M   \rt  F^*  \rt L^* \rt 0$.   
     The assertion 
of Theorem \ref{the2} follow from Theorem  \ref{e2general}.
\end{proof}

\section{First Hilbert coefficient of the canonical module}
In this section $(A,\m)$ is a \CM \ local ring with a canonical module $\om_A$.
It is well known that $e_0(\om_A) = e_0(A)$. 
 In Theorem \ref{omega} we obtain bounds for $e_1(\om_A)$.

\begin{proof}[Proof of Theorem \ref{omega}.]
Using  \cite[Theorem 18]{Pu1} it follows that $  e_1(\om_A) \leq \tau e_1(A)$ with equality
\IFF \ $A$ is Gorenstein. 
We prove the lower bound on $e_1(\om_A)$.
 There is nothing
to prove when $A$ is Gorenstein. So we assume $A$ is \emph{not}  Gorenstein.
By Remark \ref{Pr0}  we may assume $k = A/\m$
 is infinite. 

\noindent\textbf{Reduction to dimension 1:} 
Let  $x_1,\ldots,x_{d}$ be a
maximal
$\om_A\otimes A$ superficial sequence.  Set
$J = (x_1,\ldots,x_{d-1})$ and $B =  A/J$. Clearly $B$ is  a \CM \ local ring of dimension $1$.
 We also
have
 $\omega_B \cong \om_A/J\om_A$;  cf. \cite[Theorem 3.3.5]{BH}.
 By Remark \ref{Pr2}(3)
we have $e_1(B) = e_1(A)$ and $e_1(\om_B) = e_1(\om_A)$. Also 
$\type(A) = \type(B)$. 

\noindent(1.)  Set $\om = \om_A$ and 
 $M^\dagger = \Hom_A(M,\om )$.  We dualize the exact sequence
$0 \rt N \rt A^{\tau} \rt \om \rt 0$, to obtain the exact sequence
  (see \cite[Theorem 3.3.10(d)]{BH})
\begin{equation*}
0 \xar A \xar \om^{\tau} \xar N^\dagger  \xar 0. \tag{*}
\end{equation*}
 Let $x$ be $A\oplus \om \oplus N^\dagger $-superficial. For $i \geq 0$ give $L_i(N^\dagger)$,
 $L_i(\om)$, and $L_i(A)$ 
 the
 $A[X]$-module structure  as described in Remark \ref{modstruc}.

By Proposition \ref{propE}.1 we have an exact sequence
\begin{equation*}
L_1( N^\dagger) \xrightarrow{\delta} L_0(A) \rt L_0(\om^{\tau}) \rt L_0( N^\dagger ) \rt 0.\quad \text{Set}\  K = \image \delta. 
\tag{$a$}
\end{equation*}

Note that $N^\dagger$ is not free, since otherwise by (*)  we get
 $\om$ is free, a contradiction; since $A$ is \emph{not} Gorenstein.
Using Lemma \ref{dim1} there exists an $\m$-primary ideal $\q$
such that $ L_1( N^\dagger)$ is a  finitely generated  $R = A/\q[X]$- module of dimension $1$. So $K = \bigoplus_{n\geq 0}K_n$ is a finitely generated $R$-module.
Furthermore $\dim K \leq \dim R =1$. Set $\sum_{n\geq 0}\lambda(K_n)z^n = l_K(z)/(1-z)$ and $l_K(1) \geq 0$.
Note that  $l_K(1)  \neq 0$ iff $\dim K =1$.
Using $(a)$ we get 
\begin{align*}
 (1-z)l_K(z) &= - \tau h_\om(z) +  h_A(z)    + h_{ N^\dagger}(z) \quad \text{and so} \tag{$b$}    \\
\tau e_1(\om) &= e_1(A) + e_1(N^\dagger) + l_K(1).  \tag{$c$}
\end{align*}

Since all the terms involved in (c) are non-negative it follows that
\begin{equation*}
\tau e_1(\om) \geq e_1(A) \quad \text{with equality iff} \  e_1(N^\dagger) = l_K(1) = 0.  \tag{$d$}
\end{equation*}
This proves the assertion for lower bound of $e_1(\om_A)$ in (1). 

\noindent(3.)  If $G(A)$ is \CM \ then $x^*$ is $G(A)$-regular (\ref{Pr2}.5.a). Using
Proposition \ref{propE}.3 we get that $X$ is $L_0(A)$-regular. From 
$(a)$ it follows that either $K = 0$ or $X$ is an $K$-regular element.

If $K = 0$ then using $(b)$ we get that $\tau h_\om(z) =  h_A(z)    +  h_{ N^\dagger}(z)$.
So $ \tau e_i(\om_A) =  e_i(A) +  e_i(N^\dagger)$ for all $i \geq 0$. 
  By (\ref{Pr2}.6.b) $ e_i(N^\dagger) \geq 0$ for all $i \geq 0$. Thus
$ \tau e_i(\om_A) \geq  e_i(A)$ for all $i \geq 0$.

If $K \neq 0$ then $X$ is $K$-regular. Also $\dim K = 1$. Thus 
$K$ is a \CM \  $R$-module. So all 
the coefficients of $l_{K}(z)$ are non-negative. Using $(b)$ we get
\[
\tau e_i(\om_A) -  e_i(A) -  e_i(N^\dagger) = l_{K}^{(i-1)}(1)/(i-1)! \quad \geq 0
\]
for $i \geq 1$. This gives (3), since $\dim N^\dagger = 1$ and so $e_i(N^{\dagger}) \geq 0$ for each $i \geq 1$.

\

\noindent(2.) Set $k = A/\m$. If $e_1(A) = \tau e_1(\om)$ then from $(d)$ it follows that $e_1(N^\dagger)= 0$.
Therefore $N^\dagger$ (and so $N \cong (N^\dagger)^\dagger$)  are Ulrich $A$-modules. 
 So $\mu(N) = e_0(N) = (\tau -1) e_0(A)$. Set $e =  e_0(A)$. 
Let $y$ be a         $N \oplus \om \oplus A$-superficial element.
Set $(B,\n) = (A/(y), \m/(y))$. Note that $  N/xN = k^{(\tau -1) e} $ and $\om_B = \om/y\om$.
We also have 
an exact sequence 
$0 \rt k^{(\tau -1) e}  \rt B^{\tau} \rt \om_B \rt 0.$
This yields an exact sequence
\begin{equation*}
0 \xar \Hom_B(k,  k^{(\tau -1) e}) \xar \Hom_B(k,  B^{\tau})       \xar \Hom_B(k,  \om_B) \xar \cdots.   \tag{$f$}
\end{equation*}
Note that $\Hom_B(k,  B^{\tau}) \cong k^{\tau^2}$ and $ \Hom_B(k,  \om_B) \cong k$.
So by $(f)$ we get that
\[
\text{ either} \ \tau^2 = (\tau -1) e + 1 \quad \text{or} \quad \tau^2 = (\tau -1)e.
\]

\underline{Case 1:} $ \tau^2 = (\tau -1) e + 1$.

Since $\tau \neq 1$, we get $\tau =e - 1$. So $\tau = \lambda(\n)$. Thus $\socle(B) = \n$.
Therefore $\n^2 = 0$. So $B$ has minimal multiplicity and therefore $A$ also has minimal multiplicity.

\underline{Case 2:} $ \tau^2 = (\tau -1) e $.

So  $ \tau^2 - 1  = (\tau -1) e  - 1$. As  $\tau \neq 1$, we get
$\tau + 1 = e - [1/(\tau -1)]$.
Therefore $\tau \ngeq 3$. As $A$ in not Gorenstein we have
 $\tau = 2$ and $\mu(\n) \geq 2$. So $e = 4$.  There exists two possible
Hilbert series for $B$, namely

$\displaystyle{(i) \  h_B(z) = 1 + 3z \quad \text{or} \quad (ii)\  h_B(z) = 1 + 2z + z^2.}$

\noindent\underline{Claim:}  (ii) is not possible.
\textit{Proof:} Note that  $h = \mu(\m) - 1 = \mu(\n) = 2$. So 
$e = h + 2$ and $h = 2 = \tau$. Then by   \cite[Theorem 6.12]{VaSix}, we get $h_A(z) = 1 +  2z + z^3$. So $e_1(A) = 5 \neq 2e_1(\om)$, a contradiction.
Thus only $(i)$ holds and so $B$ (and therefore $A$) has minimal multiplicity.

Conversely  if
 $A$ has minimal multiplicity then $\om$ also has minimal multiplicity.  
In particular $G(A)$ and $G(\om)$ are both \CM. Say $h_A(z) = 1 + hz$. It can be checked that
$h = \tau$. It follows that $h_{\om}(z) = h + a_1z$. Since $e_0(A) = e_0(\om)$ we get
$a_1 = 1$. Therefore we get $e_1(A) = \tau e_1(\om)$.
\end{proof}

\section{generic}
\s \label{hypgen}
In this section
 $A = k[[x_1,\ldots,x_{\nu}]]/\q$ is \CM \ of dimension $d \geq 1$ and $\q \subseteq (\xbold)^2$.
 Unless stated otherwise the field $k$ is assumed to be infinite.

\s \label{m33} Let $1\leq r \leq d $.
  Set $y_j = \sum_{i =1}^{\nu}\alpha_{ij}x_i$ for 
$j = 1,\ldots r$ and $\alpha_{ij} \in k$.
We prove that for 'sufficiently general (= \sg )' $\alpha_{ij}$, the Hilbert function of
$A/(\y)$ remains the same. Notice for \sg \  $\alpha_{ij}$ we get $\y$ to be a superficial sequence.

\begin{construction}
\label{marley}
We describe a construction due to Marley \cite[p.\ 32]{MarPhd}.
 Let $y_1,\ldots,y_r$
be an $A$-superficial sequence. Let $\mathbf{K}_{\bullet}(\y)$ be the Koszul complex
and for each  $n \geq 1$ let $\K^{(n)}(\y)$ be the subcomplex
\[
  0 \rt \m^{n+1-r}K_r \rt \ldots \rt\m^nK_1 \rt \m^{n+1}K_0 \rt 0
\]
 Let $\C^{(n)}(\y)= \K(\y)/\K^{(n)}(\y)$. So $\C^{(n)} = \C^{(n)}(\y) $ is
\[
0 \rt \frac{A}{\m^{n+1-r}} \xrightarrow{\psi^{\y}_{n,r}} \ldots \rt \left(\frac{A}{\m^{n-1}}\right)^{\binom{r}{2}}\xrightarrow{\psi^{\y}_{n,2}}\left(\frac{A}{\m^{n}}\right)^{\binom{r}{1}}  \xrightarrow{\psi^{\y}_{n,1}} \frac{A}{\m^{n+1}} \rt 0. 
\]
Clearly $H_0(\C^{(n)}) = A/(\y, \m^{n+1})$. Also 
   $H_i(\C^{(n)}) = 0$ for $i \geq 1$ and $n \gg 0$ \cite[3.6]{HM}.
\end{construction}

\begin{remark}
 As $\C^{(n)}(\y)$ is a bounded complex of modules of finite length we get
\begin{equation}\label{genbruno}
\sum_{i = 0}^{r} (-1)^i \lambda\left( \C^{(n)}(\y)\right) =   \sum_{i = 0}^{r} (-1)^i \lambda\left(H_i(\C^{(n)}(\y))\right).
\end{equation}
If $f \colon \Z \rt \Z$ is a function, set $\bigtriangleup(f) = f(n) -f(n-1)$.
 Note (\ref{genbruno}) yields
\begin{align*}
\bigtriangleup^{r}\lambda\left(\frac{A}{\m^{n+1}}\right) &= \lambda\left(\frac{A}{(\y,\m^{n+1})}\right) + w(\y,n);
 \ \text{where} \\
 w(\y,n) &=  \sum_{i = 1}^{r} (-1)^i \lambda\left(H_i(\C^{(n)}(\y))\right).
\end{align*}
\end{remark}
It follows from Marley's result that $ w(\y,n) = 0$ for $n \gg 0$. An easier way to see this is by using the
Hilbert-Samuel polynomial. It can be easily checked that 
\begin{equation}
\label{post}
 w(\y,n) = 0 \ \text{ for } \ n \geq \max\{ \post(A) + r , \post(A/(\y)) \}.
\end{equation}
Recall  $\post(A)$ is the postulation number of the Hilbert-Samuel function, see (\ref{postno}).
 
\s  The  invariant $\trm(A)$ by Trivedi \cite{Trived} is convenient
 for our purposes. 
Trivedi \cite[Theorem 2]{Trived} proved that $\post(A) \leq \trm(A)$ when $\dim A \geq 0$. Also
it can be easily checked that if $y_1,\ldots,y_r$ is a superficial sequence then $\trm(A/\y) \leq \trm(A)$.

 Thus using (\ref{post}) it follows that
\begin{equation}
\label{trm}
w(\y, n) = 0 \quad \text{for} \ n \geq \trm(A) + r. 
\end{equation}

\begin{lemma}
\label{lem:suffgen}
(with hypothesis as in \ref{hypgen})
Let $1\leq r \leq d $. Then for any two sets of $ r$, \sg \ $k$-linear
combinations of $ x_1, \ldots, x_\nu $ say $y_1, \ldots, y_r$ and $
z_1, \ldots, z_r $ 

\center$\displaystyle{ H(A/(\y) ,n) = H(A/(\z),n) \quad \text{for each} \ n \geq 0}.$
\end{lemma}
\begin{proof}
 We  use \ref{m33}, \ref{marley}. We 
 prove $w(\y,n)$ is constant for \sg\ $\alpha_{ij}$.  Since $w(\y,n) = 0$ for all $n \geq \trm(A)$ for
any superficial sequence, it suffices to show for each $s \geq 1$ we have      $\dim_k H_{s}(\C^{(n)}(\y))$ is constant for \sg \ $\alpha_{ij}$. Fix an integer $s$  with $1\leq s \leq r$.
\[
\text{Consider} \quad \psi^{\y}_{n,s} \colon \left(\frac{A}{\m^{n+1 -s}}\right)^{\binom{r}{s}}   \rt \left( \frac{A}{\m^{n+2 -s}}\right)^{\binom{r}{s-1}} .
\]
\[
\text{Let}\ \  \eta_{s-1}   =  \sup\left\{\dim_k \image \psi_{n,s}^{\y} \mid  \y = y_1, \ldots, y_r
        \mbox{ is a superficial sequence} \right\}. 
\]
\begin{eqnarray*}
\text{Set} \quad y_j & = & \sum_{i = 1}^{\nu} \alpha_{ij} x_i, \quad \alpha_{ij} \in k.
\end{eqnarray*}
\textbf{Claim:} For \sg \  $\alpha_{ij}$ \ \ \ 
$\dim_k \image \psi_{ n,s}^{\y} = \eta_{s-1}.$

\s\label{done} If we prove the claim then we are done since
\begin{align*}
 \dim_k H_{s}(\C^{(n)}(\y)) &= \dim_k \ker \psi_{n,s}^{\y} - \dim_k \image \psi_{n,s+1}^{\y} \\
     &= \dim_k \left(\frac{A}{\m^{n+1 -s}}\right)^{\binom{r}{s}} - \dim_k \image \psi_{n,s}^{\y} -\dim_k \image \psi_{n,s+1}^{\y}\\
     &= \dim_k \left(\frac{A}{\m^{n+1 -s}}\right)^{\binom{r}{s}} - \eta_{s-1}- \eta_s, \quad \text{for \sg \ $\alpha_{ij}$.}
\end{align*}
\textit{Proof of Claim:} Let $\{u_1, \ldots, u_g\}$ be a $k$-basis of $A/\m^{n+1 -s}$ and let $\{v_1,
        \ldots, v_l\}$ be a $k$-basis of $A/\m^{n+2 -s}$. 
Then $u_t e_{j_1} \wedge \cdots \wedge e_{j_m}$ with  $ 1 \leq t \leq g$, and $  1 \leq
j_1 \leq \cdots \leq j_s \leq r$ is a $k$-basis of 
$\left(A/\m^{n+1-s}\right)^{\binom{r}{s}}$ and $v_t e_{j_1} \wedge \cdots \wedge
e_{j_{s - 1}}$ with $ 1 \leq t \leq l$, and $ 1 \leq
j_1 \leq \cdots \leq j_{m - 1} \leq r$ is a $k$-basis of
$\left(A/\m^{n+2-s}\right)^{\binom{r}{s-1}} $.
Set
\begin{eqnarray*}
u_t x_i & = &  \sum_{\xi = 1}^{l} a_{\xi}^{(t, i)} v_{\xi} \quad \text{for} \ t = 1, \ldots,g \ \text{and}\  i = 1,\ldots, \nu.
\end{eqnarray*}
Since 
$y_j  =  \sum_{i = 1}^{\nu}\alpha_{ij} x_i$,
we have for $t =  1, \ldots,g $ and $j =  1, \ldots, r$,
\[
u_t y_j = \sum_{i = 1}^{\nu} \alpha_{ij} \left(\sum_{\xi = 1}^{l} a_{\xi}^{(t, i)} v_{\xi} \right) = \sum_{\xi = 1}^{l}f_{\xi}^{(t, j)}(\underline\alpha)v_{\xi}\ \text{where} \  f_{\xi}^{(t, j)}(\underline\alpha)  = 
 \sum_{i = 1}^{\nu}\alpha_{ij}a_{\xi}^{(t, i)}.
\]
\begin{align*}
\text{Then}\ \ \psi_{ n,s}^{\y}(u_t e_{j_1} \wedge &\cdots \wedge e_{j_s})  = 
                \sum_{q = 1}^{s}(-1)^{q +1}u_ty_{j_q} e_{j_1} \wedge \cdots \wedge
                \hat{e}_{j_q} \wedge \cdots \wedge e_{j_s} \\
                &= \sum_{q = 1}^{s}(-1)^{q +1}\left( \sum_{\xi = 1}^{l}
                                f_{\xi}^{(t, j_q)}(\underline\alpha) v_{\xi}\right) e_{j_1} 
                                \wedge \cdots \wedge \hat{e}_{j_q} \wedge \cdots \wedge 
                                e_{j_s} \\
                 &= \sum_{q = 1}^{s} \sum_{\xi = 1}^{l} (-1)^{q +1} 
                                f_{\xi}^{(t, j_q)}(\underline\alpha) v_{\xi} e_{j_1} 
                                \wedge \cdots \wedge \hat{e}_{j_q} \wedge \cdots \wedge 
                                e_{j_s}
\end{align*}
Thus the map $\psi_{s,n}^{\y}$ can be described by a matrix of
linear forms in $\alpha_{ij}$.  Replacing the $\alpha_{ij}$ by
variables we get a matrix $E$ of linear polynomials with coefficient
in $k$.  Then by construction $U = k^{\nu} \setminus V\left(I_{\eta}(E)\right)$ is
non-empty open subset of $k^{\nu}$.  For
$\alpha_{ij} \in U$, we get  $\dim_{k} \image \psi_{ n,s}^{\y} =
\eta_{s-1}$. As indicated in \ref{done} this finishes the proof.
\end{proof}

\textbf{Application:}
We now give  an application of Lemma \ref{lem:suffgen} to bound Hilbert coefficients of a \MCM \ module $M$ if $G(M)$ is \CM.

\begin{theorem}
\label{existR}
Let $(A,\m)$ be a equicharacteristic \CM \ local ring of dimension $d > 0$. 
Then there exists a  Artinian local ring $R$ and a one dimensional  Gorenstien local ring $T$ with the property with the property

\noindent(*) If $M$ is any \MCM \ $A$-module with $G(M)$ \CM \  then 
\begin{enumerate}[\rm 1.]
\item
$e_i(M) \leq e_i(R)\mu(M) \quad \text{and} \ \chi_i(M) \leq \chi_i(R)\mu(M)  \quad \text{for all } \ i \geq 0$  
\item
$e_i(M) \leq e_i(T)\type(M)\quad \text{and} \ \chi_i(M) \leq \chi_i(T)\type(M)  \quad \text{for all } \ i \geq 0.$
\end{enumerate}
Furthermore if $A = S/\q$ where $S= k[[x_1,\ldots,x_n]]$ with  $k$ is infinite and $\q \subseteq (\xbold)^2$ 
then one can take  $R= S/\q +( v_1,\ldots,v_{d})$ and  $T = S/I +( v_1,\ldots,v_{d-1})$
where 
the $v_{i}'s$ are \sg \ linear combination of $X_1,\ldots,X_n$ and $I = (f_1, \ldots,f_r) $ is an ideal
contained in $\q$ with $\grade \q = r $ and $f_1,\ldots f_r$ is a regular sequence.
\end{theorem}
\begin{proof}
Using Remark 1.1 we may assume that the residue field $k$ of $A$ is infinite. Also we may assume $A$ is complete.
This we do. By Cohen structure theorem  we can assume $A = S/\q$ where $S= k[[x_1,\ldots,x_n]]$
and  $\q \subseteq (\xbold)^2$. 

1. Let $R = A/( v_1,\ldots,v_{d}) =  S/\q +( v_1,\ldots,v_{d}) $, where the $v_{i}'s$ are \sg \ linear combination of $x_1,\ldots,x_n$.

We know that if $\y = y_1,\ldots,y_d$ are  \sg \ linear combination of $x_1,\ldots,x_n$,
 then $ y_1,\ldots,y_d$ is an $M$-superficial sequence and  the Hilbert function of $B= A/\y A $ is equal to Hilbert function of $R$,   
by Lemma \ref{lem:suffgen}.

Set $N = M/\y M$. By Remark \ref{geqgeq} we have $e_i(N) \leq e_i(B)\mu(N)$ for each $i \geq 0$.
Note that $\mu(N) = \mu(M)$ and since $G(M)$ is \CM \ $e_i(M) = e_i(N)$, for $i \geq 0$.  
Since the Hilbert function of $B$ is same as that of $R$ we have $e_i(B) = e_i(R)$ for each $i \geq 0$.
So we get the result.

2. First note that $M$ is an \MCM \ $Q = S/I$-module and $Q$ is Gorenstein. Then one uses Theorem 3 and proves 
this assertion along the same lines as 1. 
\end{proof}

\section*{Acknowledgments}
I thank  Prof. L.~L.~Avramov and Prof B.~Ulrich for 
 many discussions.
I also thank the referee for pertinent comments

\providecommand{\bysame}{\leavevmode\hbox to3em{\hrulefill}\thinspace}
\providecommand{\MR}{\relax\ifhmode\unskip\space\fi MR }
\providecommand{\MRhref}[2]{%
  \href{http://www.ams.org/mathscinet-getitem?mr=#1}{#2}
}
\providecommand{\href}[2]{#2}

\end{document}